\documentclass{amsart}
\usepackage{amssymb,amscd,verbatim, amsthm,amsmath,amsgen,xspace}
\usepackage{graphicx, color}

\usepackage[notcite, notref]{showkeys}

\newcommand \unb[2]{\underset{#1}{{\underbrace{#2}}}}

\newcommand \fk[1]{{{\mathfrak #1}}}
\newcommand \C[1]{{\mathcal #1}}

\newcommand \bb[1]{{\mathbb #1}}

\newcommand \wti[1]{{\widetilde {#1}}}

\newcommand \bC{{\bb C}}

\newcommand \bR{{\bb R}}

\newcommand\one{1\!\!1}

\newcommand\ie{{\it i.e.,~ }}

\newcommand\eg{{\it e.g.~ }}

\newcommand\ep{{\epsilon}}
\newcommand\eps{{\epsilon}}

\newcommand\La{{\Lambda}}

\newcommand\al{{\alpha}}

\newcommand \ve{{\check e}}

\newcommand \vf{{\check f}}
\newcommand \vh{{\check h}}

\newcommand\Hom{\operatorname{Hom}}
\newcommand\Ad{\operatorname{Ad}}

\newcommand{\Kt}{\widetilde{K}}

\def\Ad{\mathop{\hbox {Ad}}\nolimits}

\def\Hom{\mathop{\hbox {Hom}}\nolimits}
\def\im{\mathop{\hbox {Im}}\nolimits}

\def\ker{\mathop{\hbox{Ker}}\nolimits}

\def\Cas{\mathop{\hbox {Cas}}\nolimits}

\newcommand{\pf}{\begin{proof}}
\newcommand{\epf}{\end{proof}}
\newcommand{\eq}{\begin{equation}}
\newcommand{\eeq}{\end{equation}}
\newcommand{\eqn}{\begin{equation*}}
\newcommand{\eeqn}{\end{equation*}}

\newcommand{\frg}{\mathfrak{g}}
\newcommand{\frh}{\mathfrak{h}}
\newcommand{\frk}{\mathfrak{k}}
\newcommand{\frl}{\mathfrak{l}}

\newcommand{\frp}{\mathfrak{p}}
\newcommand{\frq}{\mathfrak{q}}

\newcommand{\frs}{\mathfrak{s}}
\newcommand{\frt}{\mathfrak{t}}
\newcommand{\fru}{\mathfrak{u}}

\newcommand{\frgl}{\mathfrak{gl}}

\newcommand{\frsp}{\mathfrak{sp}}

\newcommand{\bbC}{\mathbb{C}}

\newcommand{\bbR}{\mathbb{R}}

\newcommand{\bbZ}{\mathbb{Z}}
\newtheorem{theorem}[equation]{Theorem}

\newtheorem{lemma}[equation]{Lemma}
\newtheorem{proposition}[equation]{Proposition}
\newtheorem{definition}[equation]{Definition}
\newtheorem{remark}[equation]{Remark}

\numberwithin{equation}{section}

\begin{document}
\today

\bigskip
\title[Dirac cohomology of unipotent representations]{Dirac cohomology of unipotent representations of $Sp(2n,\bb R)$ and $U(p,q)$}

\author{Dan Barbasch}
      \address[D. Barbasch]{Dept. of Mathematics\\
                Cornell University\\Ithaca, NY 14850, U.S.A.}
        \email{barbasch@math.cornell.edu}
        \thanks{The first author was partially supported by NSF grants DMS-0967386, DMS-0901104 and an
NSA-AMS grant.}
\author{Pavle Pand\v zi\'c}
	\address[P. Pand\v zi\'c]{Department of Mathematics, University of Zagreb, Bijeni\v cka 30, 10000 Zagreb, Croatia}
	\email{pandzic@math.hr}
	\thanks{The second author was partially supported by a grant from Ministry of science, education and sport of Republic of Croatia.}

\maketitle

\bigskip

\section{Introduction}\label{sec:1}
In this paper we will study the problem of computing the Dirac
cohomology of the special unipotent representations of the real groups
$Sp(2n,\bb R)$ and $U(p,q).$ 

We start with some background and motivation.

Let  $G$ be the real points of a linear connected reductive  group.
Its Lie algebra will be denoted by $\mathfrak g_0$. Fix a Cartan
involution $\theta$ and write $\mathfrak g_0=\mathfrak k_0 +\mathfrak s_0$ for the
Cartan decomposition. Denote by $K$ the maximal compact subgroup of $G$ with
Lie algebra $\fk k_0.$ The complexification $\mathfrak g:=(\mathfrak g_0)_\bC,$
decomposes as $\fk g=\fk k+\fk s.$ We fix a nondegenerate invariant symmetric bilinear
form $B$ on $\frg_0$, negative definite on $\frk_0$ and positive definite on $\frp_0$, 
and such that $\frk_0$ is orthogonal to $\frp_0$ with respect to $B$. 
We denote the extension of $B$ to $\frg$ by the same letter.

The \textit{Dirac Inequality} of Parthasarathy \cite{P2} plays an
important role in representation theory. We recall the basics.

The adjoint representation of $K$ on $\fk s$ lifts to
$\Ad:\wti{K}\longrightarrow Spin(\fk s_0)$, where $\wti{K}$ is the spin double cover of $K$. The Dirac operator
is defined as
\[
D=\sum_i b_i\otimes d_i \quad\in U(\fk g)\otimes C(\fk s),
\]
where $C(\fk s)$ denotes the Clifford algebra of $\fk s$ with respect to the form $B$,
$b_i$ is a basis of $\fk s$ and $d_i$ is the dual basis with respect to $B$.
$D$ is independent of the choice of the basis $b_i$ and $K-$invariant. It satisfies
\[
D^2=-(\Cas_\fk g\otimes 1 +\|\rho_\fk g\|^2)+(\Delta(\Cas_\fk k)+\|\rho_\fk k\|^2).
\]
In this formula, due to Parthasarathy \cite{P1},
\begin{itemize}
\item[-] $\Cas_\fk g$ and $\Cas_\fk k$ are the Casimir operators for $\fk g$ and
  $\fk k$ respectively,
\item[-]  $\fk h=\fk t +\fk a$ is a fundamental $\theta$-stable Cartan subalgebra
with compatible systems of positive roots for $(\fk g,\fk h)$ and $(\fk k,\fk t)$,
 \item[-] $\rho_\fk g$ and $\rho_\fk k$ are the corresponding half sums of positive roots,
\item[-] $\Delta:\fk k\to U(\fk g)\otimes C(\fk s)$ is given by $\Delta(X)=X\otimes 1+1\otimes\alpha(X)$,
where $\alpha$ is the action map $\fk k\to\fk s\fk o(\fk s)$ followed by the usual identifications
$\fk s\fk o(\fk s)\cong \bigwedge^2(\fk s)\hookrightarrow C(\fk s)$.
\end{itemize}
If $\pi$ is a $(\frg,K)-$module, then $D$ induces an operator
\[
D=D_\pi: \pi\otimes Spin\longrightarrow \pi\otimes Spin,
\]
where $Spin$ is a spin module for $C(\fk s).$
If $\pi$ is unitary, then $\pi\otimes Spin$ admits a
$K-$invariant inner product $\langle\ ,\ \rangle$ such that
$D$ is self adjoint with respect to this inner product. It follows that $D^2\geq 0$ on $\pi\otimes Spin$.
Using the above formula for $D^2$, we find that
\[
Cas_\fk g+\|\rho_\fk g\|^2\le \Delta(Cas_{\fk k})+\|\rho_\fk k\|^2
\]
on any $K-$type $\tau$ occurring in $\pi\otimes Spin.$ Another
way of putting this is
\eq
\label{dirineq}
\|\Lambda\|^2\le \|\tau +\rho_\fk k\|^2,
\eeq
for any $\tau$ occurring in $\pi\otimes Spin,$ where $\Lambda$ is the
infinitesimal character of $\pi.$ This is the Dirac inequality mentioned above.

These ideas are generalized by Vogan \cite{V} and Huang-Pand\v zi\'c \cite{HP1} as follows.
For an arbitrary admissible $(\fk g,K)$ module $\pi,$  we define {\it Dirac cohomology} of $\pi$ as
\[
H_D(\pi)=\ker D/(\ker D\cap \im D).
\]
Then $H_D(\pi)$ is a module for $\wti{K}$. If $\pi$ is unitary, $H_D(\pi)=\ker D=\ker D^2.$

The main result about $H_D$ is the following theorem conjectured by
Vogan.
\begin{theorem}\label{t:basic}
  \cite{HP1}
Assume that $H_D(\pi)$ is nonzero, and contains an irreducible $\Kt-$module with highest weight $\tau$.
Let $\Lambda\in\fk h^*$ denote the infinitesimal character of $\pi$.
Then $x\Lambda=\tau +\rho_\fk k$ for some $x$ in the Weyl
group $W=W(\fk g,\fk h).$ More precisely, there is $x\in W$ such that
$x\Lambda\mid_\fk a=0$ and $x\Lambda\mid_\fk t=\tau +\rho_\fk k.$

Conversely, if $\pi$ is unitary and
$\tau=x\Lambda-\rho_\fk k$ is the highest weight of a $\Kt-$type occuring in $\pi\otimes Spin,$ then this
$\Kt-$type is contained in $H_D(\pi).$
\end{theorem}
Note that for unitary $\pi$, the multiplicity of $\tau$ in $H_D(\pi)$ is the same as
the multiplicity of $\tau$ in $\pi\otimes Spin$.

These results might suggest that difficulties should arise in passing
between $K-$types of $\pi$ and $\Kt-$types of $\pi\otimes Spin$. For
unitary $\pi$, the situation is however greatly simplified by the
Dirac inequality. Namely, together with (\ref{dirineq}), Theorem
\ref{t:basic} shows that the infinitesimal characters $\tau+\rho_\fk k$ of
$\Kt-$types in Dirac cohomology have minimal possible norm. This means that whenever such
$E_\tau$ appears in the tensor product of a $K-$type $E_\mu$ of $\pi$ and a $\Kt-$type $E_\gamma$
of $Spin$, it necessarily appears as the PRV component \cite{PRV}, \ie
\eq
\label{PRV}
\tau = \gamma+\mu^-\qquad\text{up to } W_\fk k,
\eeq
where $\mu^-$ denotes the lowest weight of $E_\mu$.

\medskip
Assume now that $\frg$ and $\frk$ have equal rank, \eg $\frg_0=\frs\frp(2n,\bbR)$
or $\frg_0=\fru(p,q)$, and assume (as we
may) that $\Lambda$ is $\frg$-dominant.  Then the above $x$ must belong to
\[
W^1=\{w\in W\,\big|\, w\rho_\frg \text{ is $\frk$-dominant}\}.
\]
The condition that $x\Lambda$ is regular and integral for $\Kt$ puts
further restrictions on both $x$ and $\Lambda$ which will be made precise later.
We will make use of the following decomposition of the $\Kt-$module $Spin$:
\eq
\label{spindec}
Spin=\bigoplus_{\sigma\in W^1} E_{\sigma\rho_\frg-\rho_\frk}.
\eeq  

\medskip 
We now describe the special unipotent representations following
\cite{ABV}. The dual groups to $Sp(2n,\bb R)$ are the real groups
$So(p,q)$ with $p+q=2n+1$ and $p\ge q.$ Special unipotent
representations are defined as the irreducible modules which have
maximal annihilator in $U(\fk g)$ and infinitesimal character
determined by nilpotent orbits in $^\vee\!\fk g=so(2n+1,\bC)$ as follows. 

Nilpotent orbits are in 1-1 correspondence with conjugacy classes of Lie triples
$\{\ve,\vh,\vf\}$ where $\ve$ is nilpotent and $\vh$ semisimple. A
representation is called special unipotent if its infinitesimal
character is $\vh/2$ and the primitive ideal is maximal. 

In order to have nonzero Dirac cohomology, the infinitesimal character must be
conjugate to an element which is regular for $\fk k.$ This
restricts the nilpotent orbits to the ones corresponding to partitions
\[
\begin{aligned}
  &(2n+1),\\
  &(2n-2k+1,2k-1,1),\\
  &(n,n,1).
\end{aligned}
\]
If $n$ is odd, the partition $(n,n,1)$ is the same as $(2n-2k+1,2k-1,1)$
with $2k-1=n,$ but when $n$ is even it is a separate case. The partition
$(2n+1)$ corresponds to the trivial representation; it can be considered as
a special case of $(2n-2k+1,2k-1,1)$, with $k=0$. We will mostly ignore this case.

The infinitesimal characters for $(2n-2k+1,2k-1,1)$ and $(n,n,1)$, $n$ even, can be written as

\begin{align}
&\Lambda_k=(n-k,n-k-1,\dots ,1,0,-1,\dots , -k+1)\label{infcharmain}\\
&(\frac{n-1}{2},\dots ,-\frac{n-1}{2}),\quad n \text{ even.} \label{infcharother}
\end{align}

The infinitesimal characters (\ref{infcharmain}) are singular integral and we will study the corresponding
unipotent representations in detail. The infinitesimal characters (\ref{infcharother}) have the integral 
system of type $D_n$. In this last case, the
unipotent representations can be identified as  induced modules which
are irreducible, namely $Ind_{GL(n)}^{Sp(2n,\bb R)}[det]$ and
$Ind_{GL(n)}^{Sp(2n,\bb R)}[triv]$. The computation of Dirac cohomology in this case is
straightforward, and follows from the more complicated situation for infinitesimal character $\Lambda_k$.
Therefore we leave this case to the reader.

Each special unipotent representation with Dirac cohomology 
has Wave Front Set contained in the closure of
the complex nilpotent orbit of $Sp(2n,\bbC)$ with partitions where all
sizes are less than or equal to $2.$ A detailed description is in the
next section.

In the case of $U(p,q),$ the dual Lie algebra is 
$^\vee\fk g=\frg\frl(n,\bbC)$ with $n=p+q.$ Nilpotent orbits are parametrized by Jordan blocks,
equivalent to partitions of $n.$ If the partition of $n$ is $(n_1,\dots
,n_k),$ then the corresponding $^\vee h/2$ is formed of the strings
\[
\begin{aligned}
(&\frac{n_1-1}{2},\dots ,-\frac{n_1-1}{2},\\
&\frac{n_2-1}{2},\dots ,-\frac{n_2-1}{2},\\  
&\dots ,\\
&\frac{n_k-1}{2},\dots ,-\frac{n_k-1}{2}).
\end{aligned}
\]
The coordinates of the highest weight of $\tau$ are formed of all integers
or all half integers because $\tau$ must occur in the tensor product
$\mu\otimes Spin$ with $\mu$ a $K-$type of $U(p,q),$ with highest
weight formed of integers. So in order to be conjugate to a $\tau +\rho_\fk
k$, the coordinates of $\La$ have to be formed of integers and half
integers. If both integers and half integers are present, then the
partition can only have two parts $(n_1,n_2)$ with $n_1\not\equiv\ n_2
\ (mod\ 2)$.  There is only one special unipotent representation for
$U(n_1,n_2)$  which is obtained by the derived functor construction
from the trivial representation of a $\theta-$stable parabolic
subalgebra. The parameter is in the good range of \cite{KV}, so these
representations are unitary, and irreducible. We will not consider
them, as they are treatable by the same methods as the next case, and
easier. In fact the methods apply for the more general case of inducing
from a unitary character $\bb C_\xi$ in a range similar to what follows.

When the coordinates of $\La$ are formed of integers or half-integers
only, the partition has to be $(n_1,n_2)$ with $n_1\equiv n_2\ (mod\ 2).$
In this case there are many more  unipotent representations. We fix an
infinitesimal character of the form
\[
(\xi +\frac{n_1-1}{2},\dots ,\xi-\frac{n_1-1}{2}\mid \frac{n_2-1}{2},\dots ,-\frac{n_2-1}{2})
\]
The special unipotent representations are those for which $\xi =0$. We
consider unitary representations with infinitesimal
character such that $\xi +\frac{n_1-1}{2}\ge \frac{n_2-1}{2}.$  
These are parametrized by $(p_1,q_1,p_2,q_2)$ as follows. 
Let $n_1=p_1+q_1$ and $n_2=p_2+q_2$ satisfying $p_1\ge q_2,\ q_1\ge
p_2.$ Let  $\fk q=\fk l  +\fk u$ be the $\theta-$stable parabolic subalgebra 
determined by 
$$
(\unb{p_1}{1,\dots ,1},\unb{p_2}{0,\dots ,0} \mid \unb{p_1}{1,\dots
  ,1},\unb{q_2}{0,\dots ,0}).
$$ 

The representations we consider are the  $\C R_\fk q^S(\bb
C_\xi)$, where $\bbC_\xi$ denotes the character of $L$ corresponding to
\eqn
(\unb{p_1}{\xi,\dots,\xi},\unb{p_2}{0,\dots,0},\,\big|\,\unb{q_1}{\xi,\dots,\xi},\unb{q_2}{0,\dots,0}).
\eeqn
These representations are in the good range, so unitary, and irreducible.  
The WF-set is the real nilpotent orbit with partition
$$
(\unb{p_2+q_2}{2,\dots ,2},\unb{p_1+p_2-q_1-q_2}{1,\dots ,1})
$$ 
and alternating signs on the rows of size 2, $p_2$ ending in plus, $q_2$
ending in minus.

There are other interesting unipotent representations (not necessarily
special unipotent) which are genuine for double covers 
$\wti{Sp(2n,\bb R)}$ and $\wti{U(p,q)}$.  
We plan to discuss them in future research, as well as the cases of orthogonal groups.

\section{Unipotent representations of $Sp(2n,\bbR)$}\label{sec:3} 

\subsection{} We first recall some structural facts. If $G$ is (a cover of)
$Sp(2n,\bbR)$, so that $\frg_0=\frsp(2n,\bbR)$, then 
$\frg=\frsp(2n,\bbC)$ and $\frk=\frgl(n,\bbC)$. The Cartan subalgebra
for both $\frg$ and $\frk$, $\frh=\frt$, is identified with $\bbC^n$
with standard basis $e_1,\dots,e_n$. The positive roots for $\frk$ are
$e_i-e_j$, $i<j$, 
while the noncompact positive roots $\Delta^+$ for $\frg$ are
$e_i+e_j$, $i<j$, and $2e_i$. In particular, 
\[
\rho_\frg=(n,\dots,1),\qquad \rho_\frk=(\frac{n-1}{2},\dots,-\frac{n-1}{2}),
\qquad \rho_n=(\frac{n+1}{2},\dots,\frac{n+1}{2}).
\]
(The entries of $\rho_\frg$ and $\rho_\frk$ decrease by one, while the
entries of $\rho_n=\rho_\frg-\rho_\frk$ are constant.)

The Weyl group $W_\frk=W(\frk,\frt)$ consists of permutations of the coordinates,
while $W=W(\frg,\frt)$ also contains arbitrary sign changes of the
coordinates. The fundamental chamber for $\frg$ is given by the
inequalities $x_1\geq x_2\geq\dots\geq x_n\geq 0$, while the
fundamental chamber for $\frk$ is given by $x_1\geq x_2\geq\dots\geq
x_n$. (These are the closed fundamental chambers; the open ones are
given by strict inequalities.)

The subset $W^1\subset W$ may be parameterized by
$\bbZ_2^n$. Namely, for any choice of sign changes
$\eps=(\eps_1,\dots,\eps_n)$, there is a unique permutation $\tau$ of
the coordinates such that for any $\frg$-dominant $(x_1,\dots x_n)$,
$\tau(\eps_1x_1,\dots,\eps_nx_n)$ is $\frk$-dominant. We will be
slightly imprecise and identify $\eps$ with the corresponding element
of $W^1$.

Let us now examine the necessary conditions on $\Lambda=(\Lambda_1,\dots,\Lambda_n)$
so that a $(\frg,K)-$module $X$ with infinitesimal character $\Lambda$ can have
nonzero Dirac cohomology. This will also explain where the expression (\ref{infcharmain})
comes from.

First, to ensure $\frk$-integrality of $x\Lambda-\rho_\frk$, where $x\in W$, $\Lambda$ itself must
be $\frk$-integral, \ie the numbers $\Lambda_i-\Lambda_j$ must be integers.

Second, $\Lambda$ may be singular for $\fk g$, but no nonzero coordinate of (the
dominant representative of) $\Lambda$
can occur more than twice, and the coordinate 0 can appear at
most once. The $x\in W^1$ corresponding to Dirac cohomology must
then put a minus on exactly one member of each pair of repeated coordinates.

We will study the case when $\Lambda$ is integral, and the
representations are unipotent. This implies that there are no gaps, 
\ie $\Lambda_i-\Lambda_{i+1}\leq 1$ for all $i=1,\dots,n-1$. 

More precisely, if $\Lambda$ is singular, it
is conjugate to a weight of the form (\ref{infcharmain}):
\[
\Lambda_k=(n-k,n-k-1,\dots,k,k-1,\dots,-k+1)
\]
where $k\leq [\frac{n+1}{2}]$ is a positive integer. This can be
conjugated to $\frk$-dominant weights in $2^{n-2k+1}$ ways, by putting
a minus sign on any of the first $n-2k+1$ coordinates. 

The case of $\La$ regular consists of $\La=\rho_\fk g$ only,
  and the representation is the trivial module. We will mostly ignore this obvious case.

\subsection{} The unipotent representations for the infinitesimal character
$\Lambda_k$ can all be obtained by the
theta correspondence from one dimensional characters of $O(p,q).$ Our basic references
for theta correspondence are \cite{H1}, \cite{H2}, \cite{H3} and \cite{KaVe}.
We note that some of the cases we cover were also studied in \cite{NOT}.

We first study the special case  $p+q=2k\leq n$. 
Let $\eps,\eta\in\{0,1\}$. Denote by $\bbC_\eps$ the character $det^\eps$ of $O(p)$, and by
$\bbC_\eta$ the character $det^\eta$ of $O(q)$. (If $p=0$, we require $\eps=0$, and if $q=0$, 
we require $\eta=0$.)
Let $\bbC_{\eps,\eta}$ be the character of
$O(p,q)$ with restriction to $O(p)\times O(q)$ equal to $\bbC_\eps\otimes\bbC_\eta$. The
representation $X=X(p,q;\eps,\eta)$ of $G$ is obtained by theta lifting
the character $\bbC_{\eps,\eta}$ from $O(p,q)$ to $G$. 

In the following, we describe the construction more precisely.
Let $W=W_+ + W_-$ be a complex orthogonal space
where $W_\pm$ are orthogonal nondegenerate spaces of dimenensions $p$
and $q$ respectively, and let $V=\C L +\C L^\perp$ be a symplectic space
where $\C L$ and $\C L^\perp$ are Lagrangian subspaces (of dimension
$n$) in duality. Let $\Omega$ be the 
corresponding metaplectic representation. Let $G_1\times
G_2=Sp(V)\times O(W),$ and let $G_1=K_1S_1,$ $G_2=K_2S_2$ be the
Cartan decompositions with $\fk g_1=\fk k_1+\fk s_1$ and $\fk g_2=\fk
k_2+\fk s_2$ the corresponding decompositions of the complexified Lie
algebras. The representations $X(p,q;\ep,\eta)$ are equal to
$\Hom_{G_2}[\Omega,\bC_{\ep,\eta}]\cong\Hom_{K_2}[\Omega/\fk s_2\Omega,\bC_{\ep,\eta}]$
because $\fk s_2$ acts trivially on $\bC_{\ep,\eta}$. 
The same model for $\Omega$ can
be used for the representation of the pair $K_1\times U(p,q).$ Then
$\Omega$ decomposes into a direct sum 
\begin{equation}
  \label{eq:cpcttheta}
\bigoplus V(\mu)\otimes Y(\theta(\mu))  
\end{equation}
where $V(\mu)$ is a $K_1-$type, and $Y(\theta(\mu))$ is a certain
highest weight module of $U(p,q).$ The $K_1-$types $\mu$ are of the
form  
\begin{equation}
  \label{eq:k1struct}
\begin{aligned}
&(\al_1,\dots &&,\al_a&&,0&&,\dots ,&&0&&,-\beta_b,\dots ,-\beta_1)+\\
+&(k,\dots    &&,k&&,0&&,\dots,&&0&&,-k,\dots ,-k)+\\
+&(\frac{p-q}{2},\dots    &&,\frac{p-q}{2}&&,\frac{p-q}{2}&&,\dots,&&\frac{p-q}{2}&&,\frac{p-q}{2},\dots ,\frac{p-q}{2})
\end{aligned}
\end{equation}
where $a\le p,$ $b\le q,$ $\al_1\ge\dots\ge \al_a>0$ and $\beta_1\ge
\dots \ge\beta_b>0,$ and $\theta(\mu)$, the lowest $K-$type of $Y$  
is given by the second part of
formula \ref{eq:kstruct} below.

The $K_1\times\big( U(p)\times U(q))-$module
structure of $\Omega/\fk s_2\Omega$ is given by
\begin{equation}
  \label{eq:kstruct}
  \begin{aligned}
&(\al_1,\dots &&,\al_a&&,0&&,\dots ,&&0&&,-\beta_b,\dots ,-\beta_1)+\\
+&(-q,\dots    &&,-q&&,\frac{p-q}{2}&&,\dots,&&\frac{p-q}{2}&&,+p,\dots ,+p)\\
&\otimes\\
&(\al_1,\dots ,&&\al_a&&0,\dots , 0 &&\,\big|\,&&0,\dots ,0 &&-\beta_b,\dots,-\beta_1)     
  \end{aligned}
\end{equation}
where $(\al_1,\dots ,\al_a,0,\dots ,0\,\big|\, 0,\dots ,0,-\beta_b,\dots
,-\beta_1)$ is $\theta(\mu)$, the lowest $U(p)\times U(q)-$type of
$Y(\theta(\mu)).$ This uses $k=\frac{p+q}{2}.$  

To compute $\Hom_{K_2}[\Omega/\fk s_2\Omega,\bC_{\ep,\eta}]$ as a $K_1-$module, we
use the fact that $(U(p),O(p))$ and $(U(q),O(q))$ are symmetric pairs,
and Helgason's theorem. 
The end result is that the
$K_1-$types of $X(p,q,\ep,\eta)$ are of the form
\[
(\frac{p-q}{2},\dots,\frac{p-q}{2})+(\eps'+2a_1,\dots,\eps'+2a_p,0,\dots,0,-\eta'-2b_q,\dots,-\eta'-2b_1)
\]
where $a_1\geq\dots\geq a_p\geq 0$,  $b_1\geq\dots\geq b_q\geq 0$ are
nonnegative integers, and $\eps',\eta'=\eps,\eta$ if $k$ is even,
$1-\eps,1-\eta$ if $k$ is odd. 

To keep the notation simple, we write $X(p,q,\ep,\eta)$ for the
  module induced from $\bC_{\ep',\eta'}$; its $K-$structure is
\eq
\label{ktypes}
(\frac{p-q}{2},\dots,\frac{p-q}{2})+(\al_1,\dots,\al_p,0,\dots,0,-\beta_q,\dots,-\beta_1)
\eeq
where $\al_i=\ep + 2a_i$ and $\beta_j=\eta +2b_j.$ So the $\al_i$ respectively
$\beta_j$ are integers of the same parity,  $\ep$ respectively $\eta.$
These modules are unitary because they occur in the stable range of
the dual pairs correspondence, and they match unitary representations
of $O(p,q)$ \cite{L}.
 
\subsection{} 
\label{othermodules}
For $n$ odd, there is another series of special unipotent
representations $X'(p,q,\ep,\eta)$, with $p+q=2k=n+1$. Here $(\ep,\eta)$
can be $(0,0)$, $(0,1)$ or $(1,0)$ in case $p$ and $q$ are both nonzero, 
and if $p$ or $q$ is zero, then there is just one case, $(\ep,\eta)=(0,0)$. 

The infinitesimal character of $X'(p,q,\ep,\eta)$ is 
\[
\La=(k-1,k-2,\dots ,-k+1).
\]
These modules are obtained from the dual pair correspondence $Sp(2n,\bR)\times
O(p,q)$. If $p$ and $q$ are both positive, the $K-$structure of $X'(p,q,\ep,\eta)$ 
is given by the following table:

\eq
\label{otherktypes}
\begin{aligned}
    &X'(p,q,1,0): &&(\unb{p}{\frac{p-q}{2}+1},\unb{q-1}{\frac{p-q}{2}})\ \ +&
(2a_1,\dots,2a_{p},-2b_{q-1},\dots ,-2b_1);\\
    &X'(p,q,0,0): &&({\frac{p-q}{2}},\dots ,{\frac{p-q}{2}})\ \ +&
(2a_1,\dots ,2a_{p},-2b_{q-1},\dots ,-2b_1) \text{ or }\\
 &&&&+(2a_1,\dots ,2a_{p-1},-2b_{q},\dots ,-2b_1);\\
    &X'(p,q,0,1): &&(\unb{p-1}{\frac{p-q}{2}},\unb{q}{\frac{p-q}{2}-1})\ \ +&
(2a_1,\dots ,2a_{p-1},-2b_{q},\dots ,-2b_1),
\end{aligned}
\eeq
where all $a_i$ and $b_j$ are nonnegative integers. 
In each case, the first summand is the lowest $K-$type of $X'(p,q,\eps,\eta)$.
 
If $p$ or $q$ is equal to zero, then the $K-$types are:
\eq
\label{otherktypespq0}
\begin{aligned}
X'(2k,0,0,0):\quad & (k+2a_1,\dots,k+2a_n);\\
X'(0,2k,0,0):\quad & (-k-2b_n,\dots,-k-2b_1),
\end{aligned}
\eeq
where all $a_i$ and $b_j$ are nonnegative integers. In each case, the lowest $K-$type is obtained
by setting all $a_i$ respectively $b_j$ equal to zero.

From the above formulas, one can see that the modules $X'(p+1,q-1,0,1)$ and  $X'(p,q,1,0)$, respectively
$X'(2k,0,0,0)$ and $X'(2k-1,1,1,0)$, respectively $X'(0,2k,0,0)$ and $X'(1,2k-1,0,1)$ have the same $K-$types.
In fact, more is true: the identities $X'(p+1,q-1,0,1)\cong X'(p,q,1,0)$, $X'(2k,0,0,0)\cong X'(2k-1,1,1,0)$
and $X'(0,2k,0,0)\cong X'(1,2k-1,0,1)$ hold.

The modules $X'(p,q,0,1)$ and $X'(p,q,1,0)$ are unitary because they are derived
functor modules induced from a character of $U(p',q')$ in a \textit{good range} \cite{KV}. Here $U(p',q')$ is
the Levi subgroup corresponding to a $\theta$-stable parabolic subalgebra $\frq=\frl\oplus\fru$ of $\frg$.
The modules $X'(p,q,0,0)$ are constituents in the two unitarily induced
modules $Ind_P^G[det]$ and $Ind_P^G[triv]$, with $P$ the Siegel parabolic. This follows from comparing the 
$K-$types (\ref{otherktypes}) with the $K-$types of the induced modules. Thus the modules $X'(p,q,0,0)$ are unitary as well. 

\begin{remark} {\rm The modules $X(p,q,\ep,\eta)$
  have asymptotic cycle/associated cycle equal to the nilpotent orbit
  corresponding to the partition $(2,2,\dots,2)$ and signed tableau
  corresponding to $p$ pluses and $q$ minuses. The modules
  $X'(p,q,\ep,\eta)$ have asymptotic support/associated cycle
  corresponding to the signed tableau with $p$ pluses and $q-1$ minuses,
  except in the case $(\ep,\eta)=(0,0);$ in this case it is the union
  of the nilpotent orbits with $p$ pluses and $q-1$ minuses and $p-1$
  pluses and $q$ minuses.}
\end{remark}

\section{Dirac cohomology in the $Sp(2n,\bbR)$ case}\label{sec:4}

\subsection{} Recall that the infinitesimal character of $X(p,q;\eps,\eta)$ is
\[
\Lambda=\Lambda_k=(n-k,\dots,k+1,k,\dots,-k+1),
\]
where $k=\frac{p+q}{2}$.
We will refer to the last $2k$ coordinates $k,k-1,\dots,-k+1$ as the core of $\Lambda$, and to the first
$n-2k$ coordinates $n-k,n-k-1,\dots,k+1$ as the tail of $\Lambda$. 

When considering the possible
\[
\tau=x\Lambda-\rho_\frk,\qquad x\in W,
\]
we will assume that $x$ fixes the core of $\Lambda$. 

\begin{definition} 
\label{d:sp} We call $\tau=x\La-\rho_{\fk k}$ \textit{special}, if $x$ fixes the core of $\La$.
\end{definition}

Let $w_0$ be the long Weyl group
element in $W(\fk k)=S_n.$ Note that $\tau$ is special if and only if $-w_0\tau$ is not special.

The following lemma justifies why we only need to compute multiplicities of special
$\tau$ in $H_D(X(p,q,\ep,\eta).$ 

\begin{lemma} For any $\tau=x\La-\rho_\frk$, 
  \label{l:sp}
\[
[E_\tau\ :\ H_D(X(p,q,\ep,\eta))]=[E_{-w_0\tau}\ :\ H_D(X(q,p,\eta,\ep))].
\]
\end{lemma}
\pf 
There is an automorphism of $Sp(2n,\bR)$ which is $-w_0$ on the
fundamental Cartan subalgebra. It interchanges $X(p,q,\ep,\eta)$ with
$X(q,p,\eta,\ep).$ The result follows from this. 
\epf

\subsection{} Thus we will compute the multiplicity only for the case
\[
x\Lambda=(i_1,\dots,i_u,k,\dots,-k+1,-j_v,\dots,-j_1),
\]
with nonnegative integers $u$ and $v$, $u+v=n-2k$, and for some integers
\[
n-k\geq i_1>i_2>\dots>i_u\geq k+1,\qquad n-k\geq j_1>j_2>\dots>j_v\geq k+1,
\]
such that
\[
\{i_1,\dots,i_u,j_1,\dots,j_v\}=\{k+1,k+2,\dots,n-k\}.
\]

Writing $n=2k+u+v$, we can write the $i$-th coordinate of $\rho_\frk$ as 
\eq
\label{rhoksp}
(\rho_\frk)_i=k+\frac{u+v+1-2i}{2}.
\eeq
Using this and putting $i_r'=i_r-k$ and $j_s'=j_s-k$, we get
\begin{multline}
\label{tausp}
\tau=x\Lambda-\rho_\frk=(i_1'-\frac{u+v-1}{2},\dots,i_u'+\frac{u-v-1}{2},\frac{u-v+1}{2},\dots,\frac{u-v+1}{2},\\
-j_v'+\frac{u-v+1}{2},\dots,-j_1'+\frac{u+v-1}{2}).
\end{multline}
Here the first $u$ and the last $v$ coordinates correspond to the tail of $\Lambda$, and hence we will call them 
the tail of $\tau$. Analogously, the $2k$ coordinates
in between them, which are all equal to $\frac{u-v+1}{2}$, correspond to the core of $\Lambda$, and we call them
the core of $\tau$.

Let us now consider the module $X=X(p,q;\epsilon,\eta)$ described in Section \ref{sec:3}. Here $p+q=2k$, so 
the infinitesimal character of $X$ is $\Lambda$.

As we know from Theorem \ref{t:basic}, from (\ref{PRV}), and from (\ref{spindec}), the multiplicity of $\tau$ in the 
Dirac cohomology of $X$ is the number of solutions to the equation
\[
w\tau=\sigma\rho_\frg-\rho_\frk+\mu^-
\]
where $w\in W_\frk$, $\sigma\in W^1$ and $\mu^-$ is the lowest weight of a $K-$type $\mu$ of $X$.

Using the description (\ref{ktypes}) of the $K-$types of $X$, we can
rewrite this equation as
\eq
\label{eqn}
w\tau+\rho_\frk+(\frac{q-p}{2},\dots,\frac{q-p}{2})=\sigma\rho_\frg+(-\beta_1,\dots,-\beta_q,0,\dots,0,\al_p,\dots,\al_1).
\eeq
We first turn our attention to the zeros in the above expression, which appear from $(q+1)-$st to $(n-p)=(q+u+v)-$th coordinate. Using (\ref{rhoksp}) and writing $k=\frac{p+q}{2}$, we see that the $i$-th coordinate of 
$\rho_\frk+(\frac{q-p}{2},\dots,\frac{q-p}{2})$ is
$q+\frac{u+v+1-2i}{2}$. In particular, the $(q+1)-$ to 
$(q+u+v)-$th coordinates of $\rho_\frk+(\frac{q-p}{2},\dots,\frac{q-p}{2})$ are 
\[
\frac{u+v-1}{2},\frac{u+v-3}{2},\dots,-\frac{u+v-1}{2}.
\]
We see that if we take any solution of (\ref{eqn}) and denote 
the $(q+1)-$st to the $(q+u+v)-$th  coordinates of $w\tau$ by $x_1,\dots,x_{u+v}$, then
\eq
\label{part}
(x_1,\dots,x_{u+v})+(\frac{u+v-1}{2},\frac{u+v-3}{2},\dots,-\frac{u+v-1}{2})
\eeq
is a part of $\sigma\rho_\frg$ (in fact, its $(q+1)$-st to $(q+u+v)$-th coordinates). 

The conclusion is that if we fix $w$, then a necessary condition for the existence of solutions to 
(\ref{eqn}) is that the coordinates of (\ref{part}) are strictly decreasing nonzero integers between 
$-n$ and $n$, such that no two coordinates are negatives of each other.

\begin{lemma} 
\label{lemmatail}
The above condition is equivalent to $x_1,\dots,x_{u+v}$ being the tail of $\tau$. In other words, the above condition is equivalent to 
\begin{multline*}
(x_1,\dots,x_{u+v})=\\ (i_1'-\frac{u+v-1}{2},\dots,i_u'+\frac{u-v-1}{2},-j_v'+\frac{u-v+1}{2},\dots,-j_1'+
\frac{u+v-1}{2}),
\end{multline*}
\ie to the expression (\ref{part}) being equal to $(i_1',\dots,i_u',-j_v',\dots,-j_1')$.
\end{lemma}

\pf It is clear that if $x_1,\dots,x_{u+v}$ is the tail of $\tau$, \ie (\ref{part}) is equal to 
\[
(i_1',\dots,i_u',-j_v',\dots,-j_1'),
\]
then the condition is satisfied. We have to show that this is the only possibility for the condition to hold.

Since the coordinates of (\ref{part}) are strictly decreasing, the sequence $x_1,\dots,x_{u+v}$ must be decreasing. Furthermore, these are coordinates of $\tau$, and the coordinates of $\tau$ are also decreasing.
So every $x_r$ is at most equal to the $r-$th coordinate of $\tau$
(equation (\ref{tausp})), and at least equal to $(n+1-r)$-th coordinate of $\tau$. In particular, if $1\leq r\leq u$, then
\[
x_r+\frac{u+v+1-2r}{2}\leq i_r'-\frac{u+v+1-2r}{2}+\frac{u+v+1-2r}{2}=i_r',
\]
and similarly if $u+1\leq s\leq u+v$, then
\[
x_s+\frac{u+v+1-2s}{2}\geq -j'_{u+v+1-s}.
\]
Because all coordinates of (\ref{part}) are strictly decreasing, and
the above shows that they are between $i_1'\le u+v$ and $-j_1'\ge -u-v$, they must all be
of the form $\pm t$ for some $t\in\{1,2,\dots,u+v\}$. Since the number
of coordinates is exactly $u+v$, and no two can have the same absolute
value, this implies that every
$t\in\{1,2,\dots,u+v\}$ appears as a coordinate, either with a plus or
a minus sign. 

It is now clear that if $i_1'=u+v$, then the first coordinate of (\ref{part}) equals $u+v=i_1'$, and if $j_1'$ equals $u+v$, then the last coordinate of (\ref{part}) equals $-(u+v)=-j_1'$. Otherwise, all the coordinates would be strictly between $u+v$ and $-(u+v)$, and this is impossible.

If we now throw the coordinate equal to $\pm(u+v)$ out, we end up in the exactly analogous situation and we can repeat the same argument for the biggest remaining coordinate, $\pm(u+v-1)$. Continuing like this, we eventually conclude that the expression (\ref{part}) is indeed equal to $(i_1',\dots,i_u',-j_v',\dots,-j_1')$.
\epf

\subsection{} From now on we will be assuming that the condition of Lemma \ref{lemmatail} holds. This in fact fixes $w\tau$ to be
\begin{align*}
w\tau=\left(\frac{u-v+1}{2},\dots,\frac{u-v+1}{2},i_1'-\frac{u+v-1}{2},\dots,i_u'+\frac{u-v-1}{2}\right.,\\
\left.-j_v'+\frac{u-v+1}{2},\dots,-j_1'+\frac{u+v-1}{2},\frac{u-v+1}{2},\dots,\frac{u-v+1}{2}\right),
\end{align*}
where the constant $\frac{u-v+1}{2}$ appears in the first $q$ places and in the last $p$ places.

In the following we only consider the first $q$ and the last $p$ coordinates of the equation (\ref{eqn}), and we will separate them by a bar. We will denote by $\pi$ the corresponding projection from $\bbC^n$ onto $\bbC^{q+p}$. For $\rho_\frk$, we get
\[
\pi(\rho_\frk)=\left(k+\frac{u+v-1}{2},\dots,k+\frac{u+v+1}{2}-q\,\big|\,k-\frac{u+v+1}{2}-q,\dots,-k-\frac{u+v-1}{2}\right).
\]
Adding to this $\pi(w\tau)=(\frac{u-v+1}{2},\dots,\frac{u-v+1}{2})$ and $(\frac{q-p}{2},\dots,\frac{q-p}{2})$, and remembering that $k=\frac{p+q}{2}$, we see that the left hand side of (\ref{eqn}) gives
\[
(u+q,u+q-1,\dots,u+1\,\big|\,-v,-v-1,\dots,-v-p+1).
\]
The corresponding piece of the right hand side of (\ref{eqn}) is $\pi(\sigma\rho_\frg)$, with all coordinates among 
$\pm(u+v+1),\dots,\pm n$, plus 
\[
\pi\left(\mu^--(\frac{p-q}{2},\dots ,\frac{p-q}{2})\right)=(-\beta_1,\dots,-\beta_q\,\big|\,\al_p,\dots,\al_1).
\]

To continue our analysis, we first treat separately the cases when $p$ or $q$ is equal to 0. These cases are covered by the results of \cite{HPP}, but we include them here for completeness and to illustrate our method in a relatively simple situation.  

If $q=0$, then $p=2k$ and there are no coordinates to the left of the bar. Hence
\[
\pi(\sigma\rho_\frg)=(\,\big|\, -(u+v+1),-(u+v+2),\dots,-n),
\]
so (\ref{eqn}) becomes
\[
(\,\big|\, -v,-v-1,\dots,-v-2k+1)=(\,\big|\, -(u+v+1),-(u+v+2),\dots,-n)+(\,\big|\, \al_p,\dots,\al_1).
\]
This is equivalent to $\al_i=u+1$ for all $i$. 
Remembering that $\al_i=\ep+2a_i$, where $a_i$ are integers, we see that there is a unique solution if $u\equiv \ep+1 \,(mod\ 2)$,
and that there are no solutions if $u\equiv \ep \,(mod\ 2)$. 

Similarly, if $p=0$ then $q=2k$, and (\ref{eqn}) becomes
\[
(u+2k,u+2k-1,\dots,u+1 \,\big|\,)=(n,n-1,\dots,u+v+1 \,\big|\,)+(-\beta_1,\dots,-\beta_q \,\big|\,).
\]
This is equivalent to $\beta_i=v$ for all $i$. 
Remembering that $\beta_i=\eta+2b_i$, where $b_i$ are integers, we see that there is a unique solution if $v\equiv \eta \,(mod\ 2)$,
and that there are no solutions if $v\equiv \eta+1 \,(mod\ 2)$. 

\begin{proposition}
\label{propzero} 
Let $\tau$ be given by (\ref{tausp}), with $u$ and $v$ the corresponding integers. Then the multiplicity of $\tau$ in
$H_D(X(2k,0,\ep,0))$ is one if $u\equiv \ep+1 \,(mod\ 2)$, and zero if $u\equiv \ep \,(mod\ 2)$. Furthermore,
the multiplicity of $\tau$ in
$H_D(X(0,2k,0,\eta))$ is one if $v\equiv \eta \,(mod\ 2)$, and zero if $v\equiv \eta+1 \,(mod\ 2)$. 
\end{proposition}
\pf This was proved in the discussion of the previous paragraphs.
\epf

\subsection{} In the following we assume that $p$ and $q$ are both positive. The situation is now more complicated because
$\pi(\sigma\rho_\frg)$ is no longer uniquely determined. The following lemma gives some obvious necessary conditions 
on $\sigma$ solving (\ref{eqn}).

\begin{lemma}
\label{sigma_general}
Assume that (\ref{eqn}) holds. Then:
\begin{enumerate}
\item The first $q$ coordinates of $\sigma\rho_\frg$ must be of alternating parity. The same holds for the last $p$ coordinates of $\sigma\rho_\frg$.
\item The first $q$ coordinates of $\sigma\rho_\frg$ are positive, while the last $p$ coordinates are negative.
\item Either the first coordinate of $\sigma\rho_\frg$, $(\sigma\rho_\frg)_1$,  is $n$, or the last coordinate $(\sigma\rho_\frg)_n$ is $-n$. Likewise, either $(\sigma\rho_\frg)_q=u+v+1$, or $(\sigma\rho_\frg)_{n-p+1}=-(u+v+1)$. 
\end{enumerate}
\end{lemma}
\pf
(1) follows from the fact that the same statement holds for $(u+q,u+q-1,\dots,u+1\,\big|\,-v,-v-1,\dots,-v-p+1)$, while the first $q$ and the last $p$ coordinates of $(-\beta_1,\dots,-\beta_q\,\big|\,\al_p,\dots,\al_1)$ have constant parity.

(2) follows from the fact that the same statement holds for $(u+q,u+q-1,\dots,u+1\,\big|\,-v,-v-1,\dots,-v-p+1)$, while the opposite holds for $(-\beta_1,\dots,-\beta_q\,\big|\,\al_p,\dots,\al_1)$.

(3) follows from the fact that the coordinates of $\sigma\rho_\frg$ are strictly decreasing. 
\epf

In particular, it follows that we can write 
\[
\pi(\sigma\rho_\frg)=(x_1,\dots,x_q\,\big|\, -y_p,\dots,-y_1),
\]
where $n\geq x_1>x_2>\dots >x_q\geq u+v+1$ are integers of alternating parity, $n\geq y_1>y_2>\dots >y_p\geq u+v+1$ are integers of alternating parity,
and $x_i\neq y_j$ for all $i$ and $j$. So we see that solving (\ref{eqn}) is equivalent to solving 
\[
\begin{aligned}
&(u+q,\dots &&,u+1&&\,\big|\,&&-v,\dots &&,-v-p+1&&)=\\
&(x_1,\dots &&,x_q&&\,\big|\,&&-y_p,\dots &&,-y_1&&)+\\
&(-\eta-2b_1,\dots &&,-\eta-2b_q&&\,\big|\,&&\ep+2a_p,\dots &&,\ep+2a_1&&).  
\end{aligned}
\]

The lemma below shows that (\ref{eqn}) is further equivalent to the following modulo 2 version of the above equation:
\[
\begin{aligned}
&(u+q,\dots &&,u+1&&\,\big|\,&&-v,\dots &&,-v-p+1&&)+\\
&(x_1,\dots &&,x_q&&\,\big|\,&&-y_p,\dots &&,-y_1&&)=\\
&(\eta,\dots &&,\eta&&\,\big|\,&&\ep,\dots &&,\ep&&).  
\end{aligned}
\]

\begin{lemma} For a fixed $\tau$ as in (\ref{tausp}), assume that $(x_1,\dots,x_q\,\big|\, -y_p,\dots,-y_1)$ solves the above modulo two equation. 
Then there is a unique solution of
(\ref{eqn}) such that $\pi(\sigma\rho_\frg)=(x_1,\dots,x_q\,\big|\, -y_p,\dots,-y_1)$.
\end{lemma}
\pf
It is clear that the $a_i$ and $b_j$ are uniquely determined, and they are integers since the modulo two equation is satisfied. Moreover,
since $x_1>x_2>\dots >x_q$ and since $u+q,\dots,u+1$ decrease by one, it follows that $a_1\geq\dots\geq a_p$. Likewise, since $y_1>y_2>\dots >y_p$ 
and since $-v,\dots,-v-p+1$ decrease by one, it follows that $b_1\geq\dots\geq b_q$. So we indeed obtain a solution of (\ref{eqn}).
\epf

\subsection{} We are now ready to complete the computation.
Denote by $u_\tau:=0,1$ respectively $v_\tau$ the
parity of $u$ respectively $v$ occurring in the expression of $\tau.$ Recall the notion of $\tau$ special from Definition \ref{d:sp}.

\begin{theorem}
\label{thmHD}
Let $X=X(p,q;\eps,\eta)$ where $p+q=2k$. The general formula for the
cohomology is
\[
\begin{aligned}
H_D(X)=&\sum_{\tau\ special} [E_\tau :
H_D(X(p,q,\ep,\eta))]E_\tau +\\ 
+&\sum_{\tau\ not\ special} [E_{-w_0\tau} :
H_D(X(q,p,\eta,\ep))]E_\tau.
\end{aligned}
\]
For $\tau=x\La-\rho_{\fk k}$ special as
in (\ref{tausp}), the multiplicity is as follows. (For $\tau$ not special, see Lemma \ref{l:sp}.)
\begin{enumerate}
\item $p,q$ even and positive. (For $p$ or $q$ equal to zero, see Proposition \ref{propzero}.)
\begin{description}
\item[Ia] $\eps+\eta\equiv n (mod\ 2)$ and $u_\tau\equiv \eps (mod\ 2)$. The multiplicity of $\tau$ in $H_D(X)$ is
\[
\binom{k-1}{\frac{p}{2}}=\binom{k-1}{\frac{q-2}{2}}.
\]
\item[Ib] $\eps+\eta\equiv n (mod\ 2)$ and $u_\tau\equiv \eps +1(mod\ 2).$ The multiplicity of $\tau$ in $H_D(X)$ is
\[
\binom{k-1}{\frac{p-2}{2}}=\binom{k-1}{\frac{q}{2}}.
\]
\item[IIa] $\eps+\eta\equiv n+1 (mod\ 2)$ and $u_\tau\equiv \eps +1(mod\ 2).$ The multiplicity of $\tau$ in $H_D(X)$ is
\[
\binom{k}{\frac{p}{2}}.
\]
\item[IIb] $\eps+\eta\equiv n+1 (mod\ 2)$ and $u_\tau\equiv \eps (mod\ 2).$ The multiplicity of $\tau$ in $H_D(X)$ is 0.
\end{description}
\item $p,q$ odd.
\begin{description}
\item[I] $\eps+\eta\equiv n (mod\ 2).$ The multiplicity of $\tau$ in $H_D(X)$ is 
\[
\binom{k-1}{\frac{q-1}{2}}=\binom{k-1}{\frac{p-1}{2}}.
\]
\item[II] $\eps+\eta\equiv n+1 (mod\ 2).$ The multiplicity of $\tau$ in $H_D(X)$ is $0$. In this case, $H_D(X)=0$.
\end{description}
\end{enumerate}
\end{theorem}
\pf
\[
u_\tau +q+n+\eta\equiv 0 \quad\text{ or }\quad -v_\tau-p+1-n+\ep\equiv
u_\tau +q+1+\ep\equiv 0.
\]
In the following we write $u$ for $u_\tau$. 

\begin{enumerate}
\item $p,q$ even. 
\begin{description}
\item[Ia] $u+n+\eta\equiv 0,\qquad u+\ep+1\equiv 1.$ We use notation $\one$ for the vector
$(1,\dots,1\,\big|\, 1,\dots,1)$. The parities of the coordinates of $\pi(\sigma\rho_\frg)$ are 
\[
(u+\eta)\one +(0,1,\dots ,0,1\,\big|\, 0,1,\dots ,0,1),
\]
and we conclude $x_1=n.$ Then the $q-$th coordinate gives
\begin{equation}\label{eq:sigp}
(u+v+1)\equiv n-2k+1\equiv n+1\equiv u+\eta+1,
\end{equation}
so $x_q=u+v+1.$ The pairs $(n-1,n-2)\dots (u+v+3,u+v+2)$ must each
occur on one side or the other of $\,\big|\,.$ There are $q-2$ available coordinates on the left side
of $\,\big|\,$ and $p$ available coordinates on the right side of $\,\big|\,$, and $q-2$ and 
$p$ are even numbers. We conclude that the multiplicity
is $\binom{k-1}{p/2}=\binom{k-1}{(q-2)/2}$. 
\item[Ib] $u+n+\eta\equiv 1,\qquad u+\ep+1\equiv 0.$ In this case
$y_1=n$ and $y_p=u+v+1$, and the same argument as in (Ia) implies
that the multiplicity of $\tau$ is $\binom{k-1}{(p-2)/2}=\binom{k-1}{q/2}.$
\item[IIa] $u+n+\eta\equiv u+\ep+1\equiv 0.$ The parities of the coordinates of $\pi(\sigma\rho_\frg)$ are
\[
(u+\eta)\one +(0,1,\dots ,0,1\,\big|\, 1,0\dots ,1,0),
\] 
Then either $x_1=n$ or $y_1=n.$ In the first case we must have
$x_2=n-1$ in the second $y_2=n-1.$ The pairs
\[
(n,n-1)\dots (u+v+2,u+v+1)
\]
each occur on one side or the other of $\,\big|\,.$ There are
$\binom{k}{p/2}=\binom{k}{q/2}$ possible solutions.
\item[IIb] This follows from Lemma \ref{sigma_general} (3).
\end{description}

\item $p,q$ odd.
\begin{description}
\item[Ia] $u+1+n+\eta\equiv 0,\qquad u+\ep\equiv 1.$ The parities of the coordinates of $\pi(\sigma\rho_\frg)$ are
\[
(u+1+\eta)\one +(0,1,\dots ,1,0\,\big|\, 1,0\dots ,0,1),
\]
and we conclude $x_1=n.$ Then the $q-$th coordinate gives
\[
(u+v+1)=n+1\equiv u+\eta,
\]
so $y_p=u+v+1.$ The pairs 
\[
(n-1,n-2)\dots (u+v+3,u+v+2)
\]
each must occur on one side of $\,\big|\,$ only. There are an even number of coordinates,
$p-1$, $q-1$ respectively. The multiplicity is
$\binom{k-1}{(p-1)/2}=\binom{k-1}{(q-1)/2}$.
\item[Ib] $u+1+n+\eta\equiv 1,\qquad u+\ep\equiv 0.$ In this case the parities of the coordinates of $\pi(\sigma\rho_\frg)$ are
\[
(u+1+\eta)\one +(0,1,\dots ,1,0\,\big|\, 1,0,\dots ,0,1),
\]
we conclude $y_1=n,$ and the argument is the same as in (Ia).
\item[IIa] $u+1+n+\eta\equiv u+\ep\equiv 0.$ The parities of the coordinates of $\pi(\sigma\rho_\frg)$ are
\[
(u+1+\eta)\one +(0,1,\dots ,1,0\,\big|\, 0,1\dots ,1,0),
\]
and we conclude $x_1=n$ or $y_1=n$. In the first case $x_2=n-1$, and in the
second case $y_2=n-1.$ The pairs
\[
(n,n-1)\dots (u+v+2,u+v+1)
\]
each occur on one side of $\,\big|\,$ only. But $p$ and $q$ are odd, so
there is no space. There are no solutions in this case. This case is symmetric with respect to changing 
$(p,q,\ep,\eta)$ to $(q,p,\eta,\ep)$, so it follows that $H_D(X)=0$.

\item[IIb] $u+1+n+\eta\equiv u+\ep\equiv 1.$ The multiplicity of $\tau$ is 0 by Lemma \ref{sigma_general} (3).

\end{description}
\end{enumerate}
\epf

\subsection{} Let now $n$ be odd and let $X'=X'(p,q,\ep,\eta)$ be one of the additional modules described in Subsection
\ref{othermodules}. In particular, $p+q=2k=n+1$. Recall that the infinitesimal character of $X'$ is
\[
\Lambda=(k-1,\dots,-k+1)=\rho_\frk.
\]
Clearly, the only $x\Lambda$, $x\in W$, which is dominant regular for $\frk$, is $\Lambda$ itself. Since the candidates
for highest weights in $H_D(X')$ are all of the form $\tau=x\Lambda-\rho_\frk$, it follows that the only candidate is the trivial
$\Kt-$module $E_0$. To find its multiplicity, we have to count the number of solutions to the equation
\eq
\label{sigmaX'}
\sigma\rho_\frg -\rho_\frk+\mu^-=0,
\eeq
where $\sigma\in W^1$ and $\mu^-$ is the lowest weight of a $K-$type of $X'$.

\begin{theorem}
\label{thmHD'}
Let $X'=X'(p,q;\eps,\eta)$ where $p+q=2k=n+1$. Then the Dirac cohomology of $X'$ is
\[
H_D(X')=[E_0 : H_D(X')]E_0,
\]
where $E_0$ is the trivial $\Kt-$module.
The multiplicity is as follows.
\begin{enumerate}
\item $p,q$ even and positive.
\begin{description}
\item[I] $(\eps,\eta)=(1,0)$. The multiplicity of $E_0$ in $H_D(X')$ is
\[
\binom{k-1}{\frac{p}{2}}=\binom{k-1}{\frac{q-2}{2}}.
\]
\item[II] $(\eps,\eta)=(0,0)$. The multiplicity of $E_0$ in $H_D(X')$ is
\[
\binom{k}{\frac{p}{2}}=\binom{k}{\frac{q}{2}}.
\]
\item[III] $(\eps,\eta)=(0,1)$. The multiplicity of $E_0$ in $H_D(X')$ is
\[
\binom{k-1}{\frac{q}{2}}=\binom{k-1}{\frac{p-2}{2}}.
\]
\end{description}
\item $p,q$ odd.
\begin{description}
\item[I] $(\eps,\eta)=(1,0)$ or $(0,1)$. The multiplicity of $E_0$ in $H_D(X')$ is
\[
\binom{k-1}{\frac{p-1}{2}}=\binom{k-1}{\frac{q-1}{2}}.
\]
\item[II] $(\eps,\eta)=(0,0)$. The multiplicity of $E_0$ in $H_D(X')$ is 0, so $H_D(X')=0$.
\end{description}
\end{enumerate}
\end{theorem}
\pf
Let $(\ep,\eta)=(1,0)$. Then by (\ref{otherktypes}), (\ref{sigmaX'}) can be rewritten as
\[
\sigma\rho_\frg+(-2b_1,\dots,-2b_{q-1},2a_p,\dots,2a_1)=(q-1,q-2,\dots,1,-1,-2,\dots,-p).
\]
It follows that $\sigma\rho_\frg=(x_1,\dots,x_{q-1},-y_p,\dots,-y_1)$, with $x_i$ and $y_j$ positive.
Moreover, $x_1,\dots,x_{q-1}$ decrease and have the same parity as $q-1,\dots,1$, while $-y_p,\dots,-y_1$ 
decrease and have the same parity as $1,\dots,p$.

If $p$ and $q$ are even, it follows that $x_1=n$. Now we look at the pair $(n-1,n-2)$; it can either
be $(x_2,x_3)$, or $(y_1,y_2)$ (assuming $q\geq 4$ and $p\geq 2$). In either case, if there is still enough space, we again have two 
choices for the next pair, $(n-3,n-4)$, and so on. There are $k-1$ pairs, and we have to choose $p/2$ 
of them to go to the right side of $\sigma\rho_\frg$, and the remaining $(q-2)/2$ of them go to the left side. 
This proves 1) I. 

If $p$ and $q$ are odd, then $y_1=n$, and we have to place each of the $k-1$ pairs, $(n-1,n-2),(n-3,n-4),\dots$
to the left side or to the right side. The sides are now of lengths $q-1$ respectively $p-1$. This leads to the
first half of 2) I.

The reasoning for $(\ep,\eta)=(0,1)$, \ie for 1) III and the second half of 2) I, is completely analogous. 

If $(\ep,\eta)=(0,0)$, then by (\ref{otherktypes}), (\ref{sigmaX'}) can be rewritten as
\eq
\label{sigma00}
\sigma\rho_\frg+(-2b_1,\dots,-2b_{q-1},2c,2a_{p-1},\dots,2a_1)=(q-1,\dots,1,0,-1,\dots,-p+1),
\eeq
where $c$ is either $a_p$ or $-b_q$.
It follows that $\sigma\rho_\frg=(x_1,\dots,x_{q-1},z,-y_{p-1},\dots,-y_1)$, where $x_i$ and $y_j$ are positive, and $z$ can
be positive or negative.
Moreover, $x_1,\dots,x_{q-1}$ decrease and have the same parity as $q-1,\dots,1$; $-y_{p-1},\dots,-y_1$ 
decrease and have the same parity as $1,\dots,p$; and $z$ is even, $z<x_{q-1}$ if $z>0$, and $z>-y_{p-1}$ if $z<0$.

If $p,q$ are even, then to get $z>0$ we have to choose $q/2$ pairs of the possible $k-1$ pairs to go to the left (positive) side;
that is $\binom{k-1}{q/2}$ possibilities. Likewise, to get $z<0$ we have to choose $p/2$ pairs of the possible $k-1$ pairs to go 
to the right (negative) side; that is $\binom{k-1}{p/2}=\binom{k-1}{q/2-1}$ possibilities. The total number of choices is thus
\[
\binom{k-1}{q/2}+\binom{k-1}{q/2-1}=\binom{k}{q/2}.
\]
This proves 1) II. 

If $p,q$ are odd, then we see from (\ref{sigma00}) that there is no possibility to place the largest coordinate $n$, 
as both $x_1$ and $y_1$ have to be even. This proves 2) II. 
\epf

We remark that the cases $X'(2k,0,0,0)$ and $X'(0,2k,0,0)$ which appear to be missing from the above theorem, are in fact also covered.
Namely, as we remarked in Subsection \ref{othermodules}, $X'(2k,0,0,0)=X'(2k-1,1,1,0)$, while $X'(0,2k,0,0)=X'(1,2k-1,0,1)$.

\section{Unipotent representations of $U(p,q)$}
\label{repsupq}

\subsection{} Let $G=U(p,q)$ where $p\geq q\geq 0$ are integers of the same parity. We fix the maximal compact subgroup $K$ to be 
$U(p)\times U(q)$, embedded as block-diagonal matrices. We denote by $\frg=\frk\oplus\frs$ the Cartan decomposition of the complexified Lie algebra of $G$. We use the standard coordinates for the common Cartan subalgebra 
$\frh\cong\bbC^{p+q}$ of $\frg$ and $\frk$ and its dual. We will often separate the first $p$ from the last $q$ coordinates by a bar.

With the usual choices of positive roots for $\frg$ and $\frk$, we have
\[
\rho_\frg=(\frac{p+q-1}{2},\dots,-\frac{p+q-1}{2});\quad \rho_\frk=(\frac{p-1}{2},\dots,-\frac{p-1}{2}\,\big|\,\frac{q-1}{2},\dots,-\frac{q-1}{2}).
\]
The closed fundamental chamber for $\frg$ is given by inequalities $x_1\geq\dots\geq x_{p+q}$. The closed fundamental chamber for $\frk$ is given by inequalities $x_1\geq\dots\geq x_p$ and $x_{p+1}\geq\dots\geq x_{p+q}$.

The Weyl group $W=W(\frg,\frh)$ consists of permutations of all $p+q$ coordinates, while $W_\frk=W(\frk,\frh)$ consists of permutations that permute separately the first $p$ and the last $q$ coordinates. Moreover, $W^1$ consists of $(p,q)$-shuffles, \ie permutations 
$(i_1,\dots,i_p\,\big|\, j_1,\dots,j_q)$ of $(1,\dots,p+q)$ such that $i_1<\dots<i_p$ and $j_1<\dots<j_q$. 

Suppose that $p=p_1+p_2$ and $q=q_1+q_2$, and that $p_1\geq q_2$,
$q_1\geq p_2$. (These conditions were explained in the Introduction.) Denote by $a:=\frac{p_1-q_2}{2}$ and
$b:=\frac{q_1-p_2}{2}.$ Let $\frq$ be the parabolic subalgebra 
defined by 
\[
\gamma=(\unb{p_1}{1,\dots,1},\unb{p_2}{0,\dots,0}\,\big|\,\unb{q_1}{1,\dots,1},\unb{q_2}{0,\dots,0})
\]
in the sense that the roots of the Levi component $\fk l$ are those which are
$0$ on $\gamma,$ the roots in the nilradical $\fk u$  are those which
are $>0$ on $\gamma.$
The corresponding Levi subgroup $L$ of $G$ is isomorphic to $U(p_1,q_1)\times U(p_2,q_2)$. Let $\xi$ be an integer and let
$X(p_1,q_1,\xi)$ be the $(\frg,K)-$module cohomologically induced from the character of $L$ corresponding to
\eq
\label{eq:xi}
\bar\xi=(\unb{p_1}{\xi,\dots,\xi},\unb{p_2}{0,\dots,0},\,\big|\,\unb{q_1}{\xi,\dots,\xi},\unb{q_2}{0,\dots,0})
\eeq
We are using normalized induction, so that the infinitesimal character of $X(p_1,q_1,\xi)$ is
\eq
\label{eq:infchar}
\Lambda=(\frac{p_1+q_1-1}{2}+\xi,\dots,-\frac{p_1+q_1-1}{2}+\xi,\frac{p_2+q_2-1}{2},\dots,-\frac{p_2+q_2-1}{2}).
\eeq

We assume that $\xi$ satisfies
\begin{equation}
  \label{eq:unitaryirred}
  \frac{p_1+q_1-1}{2}+\xi\ge \frac{p_2+q_2-1}{2}\qquad\text{and}\qquad -\frac{p_2+q_2-1}{2}\ge
  -\frac{p_1+q_1-1}{2}+\xi.
\end{equation}
The first of these inequalities is the good range condition which insures that the derived functor module is irreducible and
unitary (see \cite{KV}). The second condition implies that the representation has maximal annihilator with the given
infinitesimal character, \ie it is a unipotent representation. The condition is equivalent to
\[
a+b\ge \xi\ge -a-b.
\]

\subsection{} To examine the $K-$type structure of $X(p_1,q_1,\xi)$, we consider the quantity
\[
\bar\xi+ 2\rho(\fk u\cap\fk s)-\rho(\fk
u)=\left(\xi+\frac{q_2-p_2}{2},\frac{p_1-q_1}{2}\,\big|\,
  \xi+\frac{p_2-q_2}{2},\frac{q_1-p_1}{2}\right).
  \] 
If this expression is $K-$dominant, then it is the lowest
$K-$type of $X(p_1,q_1,\xi)$. This happens precisely when 
\[
\xi+\frac{q_2-p_2}{2}\ge \frac{p_1-q_1}{2}\qquad\text{and}\qquad \xi+\frac{p_2-q_2}{2}\ge\frac{q_1-p_1}{2},
\]
which is equivalent to $\xi \ge a-b$ and $\xi\ge -a +b,$ \ie to $\xi\ge |a-b|.$ 
The second case is
\[
\xi+\frac{q_2-p_2}{2}< \frac{p_1-q_1}{2}\qquad\text{and}\qquad \xi+\frac{p_2-q_2}{2}\ge\frac{q_1-p_1}{2},
\]
which is equivalent to $-a+b\leq\xi<a-b$. (This is possible only when $a>b$.)
The third case is  
\[
\xi+\frac{q_2-p_2}{2}\ge \frac{p_1-q_1}{2}\qquad\text{and}\qquad \xi+\frac{p_2-q_2}{2}<\frac{q_1-p_1}{2},
\]
which is equivalent to $a-b\leq\xi<-a+b$. (This is possible only when $a<b$.)
The fourth case is
\[
\xi+\frac{q_2-p_2}{2}< \frac{p_1-q_1}{2}\qquad\text{and}\qquad \xi+\frac{p_2-q_2}{2}<\frac{q_1-p_1}{2},
\]
which is equivalent to $\xi < a-b$ and $\xi< -a +b,$ \ie to $\xi< -|a-b|.$ 
In each of these four cases, the lowest $K-$type $\mu_0$ is given respectively by 

\begin{align}
  &\left(\unb{p_1}{\xi+\frac{q_2-p_2}{2}},\frac{p_1-q_1}{2}\,\big|\,
\unb{q_1}{\xi+\frac{p_2-q_2}{2}},\frac{q_1-p_1}{2}\right)\label{eq:lkt1}\\
  &\left(\xi+\frac{q_2-p_2}{2}\,\big|\,
 \unb{p_2}{p_2-q_2+\frac{p_1-q_1}{2}},\xi+\frac{p_2-q_2}{2},\unb{q_2}{\frac{q_1-p_1}{2}}\right)\label{eq:lkt2}\\
 &\left(\unb{q_2}{{q_2-p_2}+\frac{q_1-p_1}{2}},\xi+\frac{q_2-p_2}{2},\unb{p_2}{\frac{p_1-q_1}{2}}\,\big|\,
 \xi+\frac{p_2-q_2}{2}\right)\label{eq:lkt3}\\
&\left(\unb{q_2}{q_2-p_2+\frac{q_1-p_1}{2}},\xi+\frac{q_2-p_2}{2}\,\big|\,
\unb{p_2}{{p_2-q_2}+\frac{p_1-q_1}{2}},\xi+\frac{p_2-q_2}{2}\right)
\label{eq:lkt4}
\end{align}

The $K-$structure is provided by the analogue of Blattner's formula, 
and equals 
\begin{equation}
  \label{eq:kstructure}
\mu_0+\left(a_1,\dots ,a_{q_2},0,\dots ,0,-b_{p_2},\dots ,-b_1\,\big|\,
  b_1,\dots , b_{p_2},0,\dots ,0,-a_{q_2},\dots ,-a_1\right)  
\end{equation}
with multiplicity 1.

\subsection{} We provide a different construction of the modules for which the
 $K-$structure is more apparent. Our basic references
for theta correspondence are again \cite{H1}, \cite{H2}, \cite{H3} and \cite{KaVe}. 
We note that some of the cases we cover were also studied in \cite{NOT}.

Consider the dual pair
$G_1\times G_2=U(p,q)\times U(q_2,p_2)$ in $Sp(2(p+q)(q_2+p_2),\bbR)$,
with Lie algebras
\[
\fk g_1\times \fk g_2=\fru(p,q)\times \fru(q_2,p_2),
\]
and maximal compact subgroups
\[
K_1=U(p)\times U(q);\qquad K_2=U(q_2)\times U(p_2).
\]
The Cartan decompositions of the two Lie algebras are denoted by
\[
\fk g_i=\fk k_i+\fk s_i.
\] 
The dual pair $G_1\times G_2$ fits into the following see-saw dual pair scheme:
\[
\begin{aligned}
&U(p,q)\times U(p,q) &&U(q_2,p_2)\times
U(q_2,p_2)\\
&&&\\
&G_1=U(p,q) && G_2=U(q_2,p_2) \\
&&&\\
&U(p)\times U(q)     && U(q_2)\times U(p_2)
\end{aligned}
\] 
Let $\Omega$ denote the metaplectic representation of $Sp(2(p+q)(q_2+p_2),\bbR)$. 
We are interested in the representation of $G_1$ given by 
$\Hom_{G_2}(\Omega,\bbC_{-\xi}),$ where $\bbC_{-\xi}$ denotes the character of $G_2$ with weight
$(-\xi,\dots,-\xi)$. The representations we discussed above are equal to $\Hom_{G_2}(\Omega,\bbC_{-\xi})\otimes\bbC_{\tilde\xi},$
where $\tilde\xi$ denotes the character of $G_1$ with weight $(\xi,\dots,\xi)$.

Since $\frs_2$ acts trivially on $\bbC_{-\xi}$, we can replace $\Omega$ by $\Omega/\fk s_2\Omega$
in the above $\Hom$ space. Then every $K_2-$type has a harmonic representative. 
The correspondence between harmonic $K-$types of the pair $U(m)\times U(r,s)$
is given by
\[
\begin{aligned}
&(a_1,\dots ,a_\ell ,0,\dots ,0,-b_k,\dots,-b_1)&&+(\frac{s-r}{2},\dots
,\frac{s-r}{2},\dots,\frac{s-r}{2})\\
&\longleftrightarrow\\
&(a_1,\dots ,a_\ell,0,\dots ,0\,\big|\, 0,\dots , 0,-b_k,\dots ,-b_1)&&+(\frac{m}{2},\dots ,\frac{m}{2}\,\big|\, -\frac{m}{2},\dots,-\frac{m}{2}).
\end{aligned}
\]
So the $K_1\times (K_2\times K_2)-$structure of $\Omega/\fk s_2\Omega$ is given by
\begin{equation}
  \label{eq:kstructure1}
  \begin{aligned}
&\left(\mu_++\frac{q_2-p_2}{2},\frac{q_2-p_2}{2},\mu_-+\frac{q_2-p_2}{2}\right)\otimes
 \left(\eta_++\frac{p_2-q_2}{2},\frac{p_2-q_2}{2},\eta_-+\frac{p_2-q_2}{2}\right)\\
&\otimes\\
&\left(\mu_+ +\frac{p_1+p_2}{2}\ \left|\
    \mu_--\frac{p_1+p_2}{2}\right)\otimes
  \left(\eta_--\frac{q_1+q_2}{2}\ \right|\ \eta_++\frac{q_1+q_2}{2}\right).     
  \end{aligned}
\end{equation}
Here $\mu_+,\eta_+$ must have nonnegative entries only, while $\mu_-,\eta_-$
must have nonpositive entries only.
The only $K_1-$types which contribute to $(-\xi\,\big|\,-\xi)$ are those
for which
\begin{align}
&\Hom_{K_2}\left[\mu_+\otimes \eta_- :-\xi
-\frac{p_1+p_2}{2}+\frac{q_1+q_2}{2}\right]\ne 0,\label{eq:kcondition1}\\
&\Hom_{K_2}\left[\mu_-\otimes \eta_+ :-\xi
+\frac{p_1+p_2}{2}-\frac{q_1+q_2}{2}\right]\ne 0, \label{eq:kcondition2} 
\end{align}
The right hand sides are $-\xi-a +b$ and $-\xi+a-b.$ The four cases of
lowest $K-$types are given by setting one of $\mu_+,\eta_-$ and one of
$\mu_-,\eta_+$ equal to zero, and the other to the right hand side,
depending on the signs of $\xi\pm a\mp b$. The $K-$structure is
determined by the equations  (\ref{eq:kcondition1}) and (\ref{eq:kcondition2}).

\section{Dirac cohomology in the $U(p,q)$ case}

\subsection{} Recall that to find the Dirac cohomology of $X(p_1,q_1,\xi)$ we first have to determine all $x\Lambda$, $x\in W$, which are dominant and regular for $\frk$. Then for each such $x\Lambda$, the $\tilde K-$type with highest weight 
$\tau=x\Lambda-\rho_\frk$ appears in $H_D(X(p_1,q_1,\xi))$ with multiplicity equal to the number of solutions $(w\tau,\mu,\sigma)$
of the equation
\eq
\label{dcoheq}
w\tau=\mu^- + \sigma\rho_\frg-\rho_\frk,
\eeq
where $w\in W_\frk$, $\mu^-$ is the lowest weight of a $K-$type of $X(p_1,q_1,\xi)$, and $\sigma\in W^1$.

In our case, the possible $x\Lambda$ are
\begin{multline}
\label{xlambda}
x\Lambda=(i_1,\dots,i_r,\frac{p_2+q_2-1}{2},\dots,-\frac{p_2+q_2-1}{2},j_1,\dots,j_s\,\big|\, \\
k_1,\dots,k_t,\frac{p_2+q_2-1}{2},\dots,-\frac{p_2+q_2-1}{2},l_1,\dots,l_u),
\end{multline}
where $i_1>\dots >i_r;\,k_1>\dots>k_t$ form a shuffle of $\frac{p_1+q_1-1}{2}+\xi,\dots,\frac{p_2+q_2-1}{2}+1$,
and $j_1>\dots >j_s;\,l_1>\dots>l_u$ form a shuffle of $-\frac{p_2+q_2-1}{2}-1,\dots,-\frac{p_1+q_1-1}{2}+\xi$.
Here $r,s,t$ and $u$ are integers satisfying 
\begin{eqnarray}
\label{rstu}
r+t=&\frac{p_1+q_1-1}{2}+\xi-\frac{p_2+q_2-1}{2} \\
s+u=&-\frac{p_2+q_2-1}{2}+\frac{p_1+q_1-1}{2}-\xi \notag \\
r+s=&p-p_2-q_2=p_1-q_2 \notag \\
t+u=&q-p_2-q_2=q_1-p_2. \notag 
\end{eqnarray}

To write down the corresponding $\tau$, it will be convenient to write 
\eq
\label{rhokupq}
\rho_\frk=(\frac{p+1}{2},\dots,\frac{p+1}{2}\,\big|\,\frac{q+1}{2},\dots,\frac{q+1}{2})-(1,2,\dots,p\,\big|\,1,2,\dots,q).
\eeq
Thus $\tau=\tilde\tau-(\frac{p+1}{2},\dots,\frac{p+1}{2}\,\big|\,\frac{q+1}{2},\dots,\frac{q+1}{2})$, where
\begin{multline}
\label{tauupq}
\tilde\tau= (i_1+1,\dots,i_r+r,\frac{p_2+q_2+1}{2}+r,\dots,\frac{p_2+q_2+1}{2}+r,j_1+p-s+1,\dots,j_s+p\,\big|\,\\
             k_1+1,\dots,k_t+t,\frac{p_2+q_2+1}{2}+t,\dots,\frac{p_2+q_2+1}{2}+t,l_1+q-u+1,\dots,l_u+q).
\end{multline}
We call the two constant strings in the above expression the core of $\tau$, and the rest of the expression the tail 
of $\tau$.

\subsection{} The $K-$types of $X(p_1,q_1,\xi)$ are given by (\ref{eq:kstructure}), with $\mu_0$ equal to one of (\ref{eq:lkt1}) - (\ref{eq:lkt4}). 
We write things out in detail for $\mu_0$ given by (\ref{eq:lkt1}), and comment what happens in other cases. 

If $\mu_0$ is given by (\ref{eq:lkt1}), then 
\begin{multline*}
\mu=(\xi+\frac{q_2-p_2}{2}+a_1,\dots,\xi+\frac{q_2-p_2}{2}+a_{q_2},
\unb{p_1-q_2}{\xi+\frac{q_2-p_2}{2},\dots,\xi+\frac{q_2-p_2}{2}},\\
\frac{p_1-q_1}{2}-b_{p_2},\dots,\frac{p_1-q_1}{2}-b_1\,\big|\,\xi+\frac{p_2-q_2}{2}+b_1,\dots,
\xi+\frac{p_2-q_2}{2}+b_{p_2},\\
\unb{q_1-p_2}{\xi+\frac{p_2-q_2}{2},\dots,\xi+\frac{p_2-q_2}{2}},\frac{q_1-p_1}{2}-a_{q_2},\dots,\frac{q_1-p_1}{2}-a_1).
\end{multline*}
Using the expression (\ref{rhokupq}) for $\rho_\frk$, we see that $\mu^--\rho_\frk=\tilde\mu-(\frac{p+1}{2},\dots,\frac{p+1}{2}\,\big|\,\frac{q+1}{2},\dots,\frac{q+1}{2})$, where
\begin{multline}
\label{mu}
\tilde\mu=(\frac{p_1-q_1}{2}-b_1+1,\dots,\frac{p_1-q_1}{2}-b_{p_2}+p_2,\xi+\frac{q_2-p_2}{2}+p_2+1,\dots,\xi+\frac{q_2-p_2}{2}+p-q_2,\\
\xi+\frac{q_2-p_2}{2}+a_{q_2}+p-q_2+1,\dots,\xi+\frac{q_2-p_2}{2}+a_1+p\,\big|\,\frac{q_1-p_1}{2}-a_1+1,\dots,\frac{q_1-p_1}{2}-a_{q_2}+q_2,\\
\xi+\frac{p_2-q_2}{2}+q_2+1,\dots,\xi+\frac{p_2-q_2}{2}+q-p_2,\xi+\frac{p_2-q_2}{2}+b_{p_2}+q-p_2+1,\dots,\xi+\frac{p_2-q_2}{2}+b_1+q)
\end{multline}

Since $(\frac{p+1}{2},\dots,\frac{p+1}{2}\,\big|\,\frac{q+1}{2},\dots,\frac{q+1}{2})$ is invariant for $W_\frk$, the equation
(\ref{dcoheq}) is equivalent to the equation
\eq
\label{dcoheqnor}
w\tilde\tau=\tilde\mu+\sigma\rho_\frg,
\eeq
with $w,\sigma$ and $\mu$ as before. Since the components $\tilde\mu_{p_2+1},...,\tilde\mu_{p-q_2}$ are increasing by 1, and since the corresponding components of $\sigma\rho_\frg$ are decreasing by at least 1, we conclude that
\eq
\label{tail1}
(w\tilde\tau)_{p_2+1}\geq\dots\geq(w\tilde\tau)_{p-q_2}.
\eeq
For the same reasons, 
\eq
\label{tail2}
(w\tilde\tau)'_{q_2+1}\geq\dots\geq(w\tilde\tau)'_{q-p_2}.
\eeq
Here and in the following we are using notation $\nu=(\nu_1,\dots,\nu_p\,\big|\,\nu'_1,\dots,\nu'_q)$ for any 
$\nu\in\bbC^{p+q}$.
\begin{lemma}
\label{tailtau}
$(w\tilde\tau)_{p_2+1},\dots,(w\tilde\tau)_{p-q_2}$ and $(w\tilde\tau)'_{q_2+1},\dots,(w\tilde\tau)'_{q-p_2}$
are exactly the components of the tail of $\tau$.
\end{lemma}
\pf
It follows from (\ref{tail1}) and (\ref{tail2}) that
\begin{eqnarray*}
(w\tilde\tau)_{p_2+1}\leq \tilde\tau_1=i_1+1,&(w\tilde\tau)_{p_2+2}\leq \tilde\tau_2=i_2+2,\dots&\qquad\text{and} \\
(w\tilde\tau)'_{q_2+1}\leq \tilde\tau'_1=k_1+1,&(w\tilde\tau)'_{q_2+1}\leq \tilde\tau'_2=k_2+2,\dots&
\end{eqnarray*}
Now either $i_1$ or $k_1$ is equal to the biggest component of $\Lambda$, $\frac{p_1+q_1-1}{2}+\xi$, while the other is strictly smaller. Assume that $i_1=\frac{p_1+q_1-1}{2}+\xi$, the other case being analogous.

It now follows from (\ref{dcoheqnor}) and (\ref{mu}) that for the corresponding components of $\sigma\rho_\frg$ we have
\begin{eqnarray*}
(\sigma\rho_\frg)_{p_2+1}\leq i_1+1-\xi-\frac{q_2-p_2}{2}-p_2-1=\frac{p_1+q_1-q_2-p_2-1}{2} &\quad\text{and}\\
(\sigma\rho_\frg)'_{q_2+1}\leq k_1+1-\xi-\frac{p_2-q_2}{2}-q_2-1<\frac{p_1+q_1-q_2-p_2-1}{2}.
\end{eqnarray*}
However, $(\rho_\frg)_{p_2+q_2+1}=\frac{p+q+1}{2}-p_2-q_2-1$ is exactly equal to $\frac{p_1+q_1-q_2-p_2-1}{2}$. It therefore follows that $(\sigma\rho_\frg)_1,\dots,(\sigma\rho_\frg)_{p_2}$ and $(\sigma\rho_\frg)'_1,\dots,(\sigma\rho_\frg)'_{q_2}$ must exactly exhaust all the components of $\rho_\frg$ that are bigger than $\frac{p_1+q_1-q_2-p_2-1}{2}$, while $(\sigma\rho_\frg)_{p_2+1}$
must be equal to $\frac{p_1+q_1-q_2-p_2-1}{2}$. This in turn means that $(w\tilde\tau)_{p_2+1}=i_1+1$.

We can repeat the same argument for the next largest component of $w\tilde\tau$ and so on, until we conclude that  
\begin{eqnarray*}
(w\tilde\tau)_{p_2+1}=i_1+1,\dots,(w\tilde\tau)_{p_2+r}=i_r+r,& \qquad\text{and}   \\ 
(w\tilde\tau)_{q_2+1}=k_1+1,\dots,(w\tilde\tau)_{q_2+t}=k_t+t.&
\end{eqnarray*}
Then we can use a similar argument bounding the components of $w\tilde\tau$ from below to conclude the rest of the lemma.
\end{proof} 

\begin{remark} {\rm The conclusion of Lemma \ref{tailtau} remains valid, with the same proof, if $\mu_0$ is given by 
(\ref{eq:lkt2}) - (\ref{eq:lkt4}). Namely, the above arguments depended only on the $(p_2+1)$-st to the $(p-q_2)$-nd 
components of the left side of $\tilde\mu$, and on the $(q_2+1)$-st to the $(q-p_2)$-nd components of the right side 
of $\tilde\mu$. These components are however the same in all four cases.}
\end{remark}

\subsection{} So we see that equation (\ref{dcoheqnor}) can have solutions only for 
\begin{multline*}
w\tilde\tau=(\unb{p_2}{\frac{p_2+q_2+1}{2}+r},i_1+1,\dots,i_r+r,j_1+p-s+1,\dots,j_s+p,\unb{q_2}{\frac{p_2+q_2+1}{2}+r}\,\big|\,\\ \unb{q_2}{\frac{p_2+q_2+1}{2}+t},k_1+1,\dots,k_t+t,l_1+q-u+1,\dots,l_u+q,\unb{p_2}{\frac{p_2+q_2+1}{2}+t}.
\end{multline*}
It remains to count the possibilities for $\sigma\rho_\frg$ and $\tilde\mu$. The result is

\begin{theorem}
\label{multtau} Let $X(p_1,q_1,\xi)$ be any of the representations considered in Section \ref{repsupq}. Then the Dirac cohomology
of $X(p_1,q_1,\xi)$ consists of $\Kt-$modules $E_\tau$ for $\tau$ defined as in (\ref{tauupq}), each appearing with multiplicity 
$\binom{p_2+q_2}{p_2}$.
\end{theorem}
\pf
We write out the proof in case $\mu_0$ is given by (\ref{eq:lkt1}). All other cases are completely analogous and left to the reader.

As we have seen in the proof of Lemma \ref{tailtau}, the largest $p_2+q_2$ components of $\rho_\frg$ are distributed between
$(\sigma\rho_\frg)_1,\dots,(\sigma\rho_\frg)_{p_2}$ and $(\sigma\rho_\frg)'_1,\dots,(\sigma\rho_\frg)'_{q_2}$, while the smallest 
$p_2+q_2$ components of $\rho_\frg$ are distributed between
$(\sigma\rho_\frg)_{p-q_2+1},\dots,(\sigma\rho_\frg)_p$ and $(\sigma\rho_\frg)'_{q-p_2+1},\dots,(\sigma\rho_\frg)'_q$.

We claim that any possible choice of strings $(\sigma\rho_\frg)_1,\dots,(\sigma\rho_\frg)_{p_2}$ and 
$(\sigma\rho_\frg)'_1,\dots,(\sigma\rho_\frg)'_{q_2}$, that is, any $(p_2,q_2)$-shuffle, leads to a unique solution of 
(\ref{dcoheqnor}). This clearly implies the statement of the theorem.

To prove this claim, let us fix the above strings in $\sigma\rho_\frg$. For any $i$ between 1 and $p_2$, the $i$th component of the equation (\ref{dcoheqnor}) reads
\eq
\label{bs}
\frac{p_1-q_1}{2}-b_i+i+(\sigma\rho_\frg)_i=\frac{p_2+q_2+1}{2}+r.
\eeq
This determines $b_i$. On the other hand, the $(p+q-i+1)-$st component of the equation (\ref{dcoheqnor}) reads
\[
\xi+\frac{p_2-q_2}{2}+b_i+q-i+1+(\sigma\rho_\frg)'_{q-i+1}=\frac{p_2+q_2+1}{2}+t.
\]
Using (\ref{rstu}), we see after some simplification that these two equations are equivalent when 
$(\sigma\rho_\frg)'_{q-i+1}=-(\sigma\rho_\frg)_i$ (and impossible otherwise). A completely analogous argument shows
that knowing $(\sigma\rho_\frg)'_j$, $j=1,\dots,q_2$, determines the $a_j$ and forces 
$(\sigma\rho_\frg)_{p-j+1}=-(\sigma\rho_\frg)'_j$.

It only remains to see that the $b_i$ and $a_j$ obtained from $(\sigma\rho_\frg)_i$ respectively $(\sigma\rho_\frg)'_j$
are always in descending order. For the $b_i$, this follows readily from (\ref{bs}) and the fact that the
$(\sigma\rho_\frg)_i$ are strictly decreasing. For the $a_j$, the argument is analogous.
\end{proof}

\end{document}